\documentclass[12pt]{amsart}
\usepackage{psfig}

\newcommand{\C}{{\mathbb C}}
\newcommand{\Q}{{\mathbb Q}}
\newcommand{\Qbar}{\overline{\Q}}
\newcommand{\Z}{{\mathbb Z}}
\newcommand{\R}{{\mathbb R}}
\newcommand{\F}{{\mathbb F}}
\newcommand{\PP}{{\mathbb P}}
\newcommand{\tors}{_{\text{tors}}}
\newcommand{\PrePer}{\operatorname{PrePer}}
\newcommand{\End}{\operatorname{End}}
\newcommand{\Gal}{\operatorname{Gal}}

\newcommand{\Fq}{{{\mathbb F}_q}}
\newcommand{\defequal}{\stackrel{\text{def}}{=}}

\newcommand{\isom}{\cong}
\newcommand{\calC}{{\mathcal C}}
\newcommand{\calM}{{\mathcal M}}
\newcommand{\calO}{{\mathcal O}}
\newcommand{\tildeD}{{\widetilde D}}

\newtheorem{cor}{Corollary}
\newtheorem{prop}{Proposition}
\newtheorem{theorem}{Theorem}

\theoremstyle{definition}
	
\newtheorem{conj}{Conjecture}

\theoremstyle{remark}

\begin{document}

\title[Rational Preperiodic Points]{The Complete Classification\\
of Rational Preperiodic Points\\
of Quadratic Polynomials over $\Q$:\\A Refined Conjecture}
\author{Bjorn Poonen}
\thanks{This research was supported by an NSF Mathematical Sciences Postdoctoral Research Fellowship.  Most of it was done at MSRI, where research is supported in part by NSF grant DMS-9022140.}
\address{Department of Mathematics, Princeton University, Princeton, NJ 08544-1000, USA}
\email{poonen@@math.princeton.edu}
\date{November 6, 1995}

\begin{abstract}
We classify the graphs that can occur as the graph of rational preperiodic points of a quadratic polynomial over $\Q$, assuming the conjecture that it is impossible to have rational points of period $4$ or higher.
In particular, we show under this assumption that the number of preperiodic points is at most~$9$.
Elliptic curves of small conductor and the genus~$2$ modular curves $X_1(13)$, $X_1(16)$, and $X_1(18)$ all arise as curves classifying quadratic polynomials with various combinations of preperiodic points.
To complete the classification, we compute the rational points on a non-modular genus~$2$ curve by performing a $2$-descent on its Jacobian and afterwards applying a variant of the method of Chabauty and Coleman.
\end{abstract}

\maketitle

\section{Introduction}
\label{intro}

Let $f:\PP^n \rightarrow \PP^n$ be a morphism of degree $d \ge 2$ defined over a number field $K$.
A point $P \in \PP^n(K)$ is called {\em periodic} (resp. {\em preperiodic}) if the sequence
	$$P,f(P),f(f(P)),f(f(f(P))),\ldots$$
is periodic (resp. eventually periodic).
The set $\PrePer(f,K)$ of preperiodic points of $f$ defined over $K$ can be made a directed graph by drawing an arrow from $P$ to $f(P)$ for each $P$.
Northcott used height functions to prove that $\PrePer(f,K)$ is always finite.
Moreover, it can be computed effectively given $f$.
These facts have been rediscovered (in varying degrees of generality) by many authors~\cite{narkiewicz}, \cite{lewis}, \cite{callsilverman}.

It is much more difficult to obtain uniform results for morphisms of fixed degree.
Morton and Silverman~\cite{mortonsilverman} have proposed the following conjecture.
\begin{conj}
\label{bigconjecture}
There exists a bound $B=B(D,n,d)$ such that if $K/\Q$ is a number field of degree $D$, and $f : \PP^n \rightarrow \PP^n$ is a morphism of degree $d \ge 2$ defined over $K$, then $\#\PrePer(f,K) \le B$.
\end{conj}
As pointed out by Silverman in talks on the subject, the special case $D=1$, $n=1$, $d=4$ is enough to imply the recently proven ``strong uniform boundedness conjecture''~\cite{merel}, since torsion points of elliptic curves are exactly the preperiodic points of the multiplication-by-2 map, and their $x$-coordinates are preperiodic points for the degree~$4$ rational map that gives $x(2P)$ in terms of $x(P)$.
A similar conjecture for polynomials over $\Fq(T)$ and its finite extensions would imply the uniform boundedness conjecture for Drinfeld modules~\cite{poonen}, which is still open.

Thus it should come as no surprise that even the simplest cases of the conjecture seem to be difficult.
For the case of quadratic polynomials over $\Q$, the following conjecture has been made~\cite{fivecycle}:
\begin{conj}
\label{periodissmall}
If $N \ge 4$, then there is no quadratic polynomial $f(z) \in \Q[z]$ with a rational point of exact period~$N$.
\end{conj}
Conjecture~\ref{periodissmall} has been verified for $N=4$ and $N=5$ (see~\cite{morton4} and~\cite{fivecycle}, respectively), and~\cite{fivecycle} presents some evidence that it holds for $N=6$ as well.

The purpose of this paper is to refine the conjecture by listing the directed graphs that occur as $\PrePer(f,\Q)$ for a quadratic polynomial $f(z) \in \Q[z]$, assuming that Conjecture~\ref{periodissmall} holds.
This list appears in Figure~1, which we will explain further in Section~\ref{preperiodic}.
It can be thought of as the analogue of the list of possible torsion sugroups of elliptic curves over $\Q$ conjectured by Levi~\cite{levi} and later by Ogg~\cite{ogg}, and finally proved by Mazur~\cite{mazur}.
In particular, we will show that Conjecture~\ref{periodissmall} implies $\#\PrePer(f,\Q) \le 9$.

\begin{figure}[p]
\centerline{\psfig{file=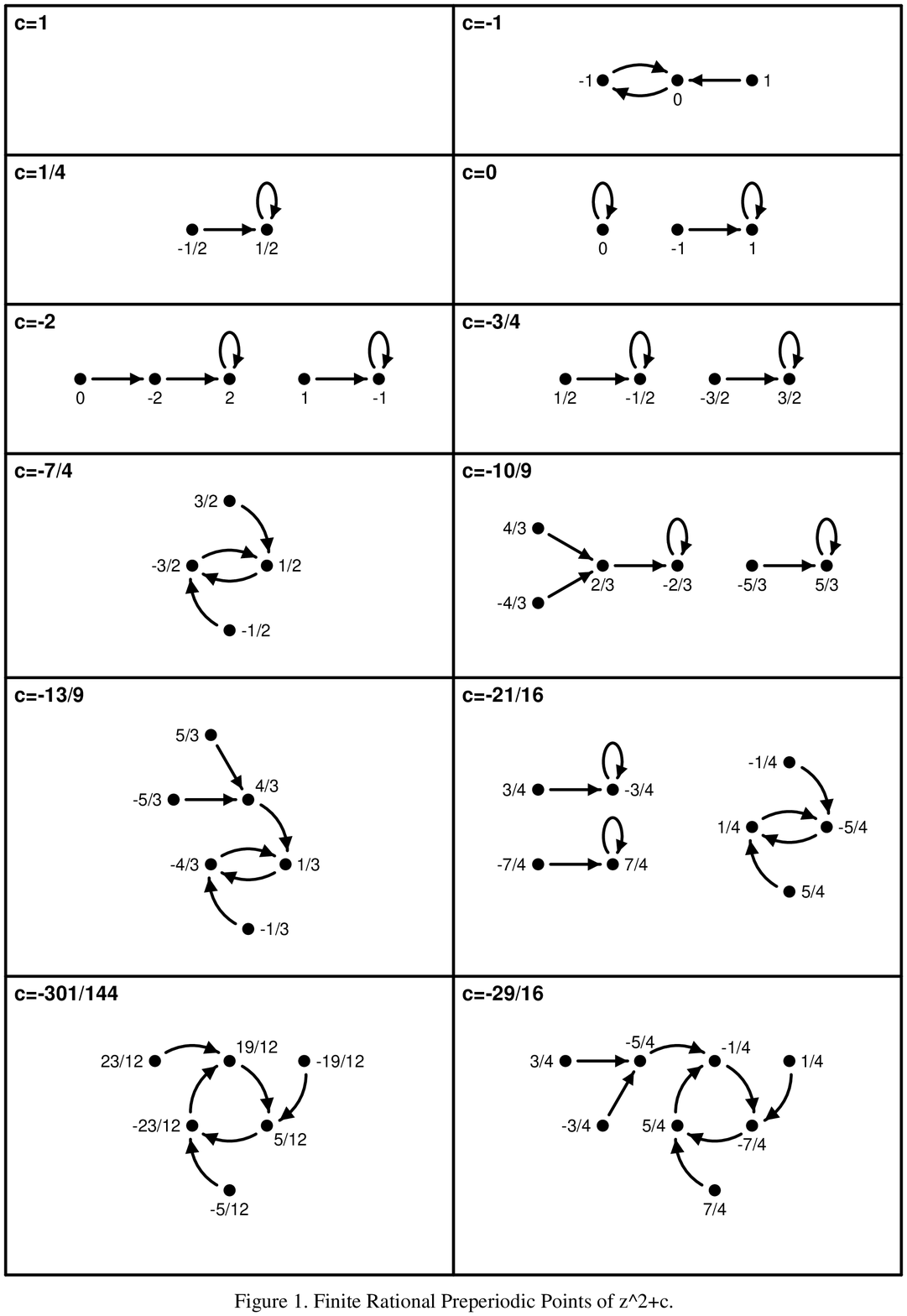,height=8in}}
\end{figure}

One can decide whether a given graph can occur in much the same way that one decides whether a given group occurs as a torsion subgroup of an elliptic curve $E$ over $\Q$: the pairs $(f(z),G)$ where $f(z)$ is a quadratic polynomial (up to linear conjugacy) and $G$ is a graph of preperiodic points (of specified shape) correspond to the points on an algebraic curve, just as elliptic curves together with level structure are parameterized by points on modular curves.
For example, quadratic polynomials together with a point of period $N$ are classified by a curve $C_1(N)$ (see~\cite{fivecycle}), which may be considered the analogue of $X_1(N)$.
The curve $C_1(N)$ also has a quotient $C_0(N)$ (the analogue of $X_0(N)$), whose points correspond to quadratic polynomials with a (Galois-stable) $N$-cycle.
To decide whether a given graph occurs, it suffices to find the rational points on the corresponding curve.

There exist also some subtle (and surprising) connections between preperiodic points of quadratic polynomials and torsion points on elliptic curves.
Morton~\cite{morton4} noticed that $C_1(4)$ was birational to $X_1(16)$, and we will show below that the curve that classifies quadratic polynomials with points of period~$2$ and~$3$ both is birational to $X_1(13)$!
Similarly, the curve that classifies quadratic polynomials with points of period~$1$ and~$3$ is birational to $X_1(18)$.
This is not a general phenomenon, however: $C_1(5)$ and $C_0(5)$ are not modular~\cite{fivecycle}.

The curves whose rational points we will need to determine will all be of genus $0$, $1$, or~$2$.
The genus~$0$ and~$1$ curves we encounter always have rational points ``at infinity,'' so each is $\PP^1$ or an elliptic curve over $\Q$.
The conductor of each elliptic curve is small (never more than $40$) and the rank is always zero (luckily), so there is no difficulty in listing the rational points in these cases.
As for the genus~$2$ curves, all except one are recognized as modular curves, whose rational points have been previously computed.
The only difficult case is the curve $C_1(3_2)$ that classifies quadratic polynomials together with a point that in two steps enters a $3$-cycle.
The curve $C_1(3_2)$ is not modular, it has no rational Weierstrass points, and its Jacobian $J$ is absolutely simple, all of which rule out possible shortcuts to finding the rational points.
We will need first to perform a general $2$-descent on $J$ to find its Mordell-Weil rank.
The rank turns out to be~$1$, so to find the rational points on $C_1(3_2)$ we apply the method of Chabauty and Coleman~(\cite{CHA},~\cite{COL}).
The refinements of the method introduced in~\cite{fivecycle} are sufficient for obtaining a tight bound on the number of rational points in most of the mod~$3$ residue classes, but for one residue class (a Weierstrass point mod~$3$), those techniques would require more $3$-adic precision than is easy to compute.
We circumvent the problem by introducing a variant of the method.

Although it would certainly be possible to do similar calculations for other number fields $K$, the results obtained might not be as conclusive.
Problems might arise in the application of the method of Chabauty and Coleman when computing the $K$-rational points of curves of genus $2$ or more, because the Mordell-Weil rank may exceed the genus.

\section{Periodic Points}
\label{periodic}

Polynomials $f(z),g(z) \in \Q[z]$ are {\em linearly conjugate} over $\Q$ if there exists a linear polynomial $\ell(z)$ such that $\ell(f(\ell^{-1}(z)))=g(z)$.
In this case, $\ell$ maps the rational preperiodic points of $f(z)$ bijectively to the rational preperiodic points of $g(z)$, also preserving the graph structure.
Every quadratic polynomial in $\Q[z]$ is linearly conjugate over $\Q$ to one of the form $z^2+c$ with $c \in \Q$, so from now on, we will only consider $f(z)=z^2+c$.
In the subsequent theorems for polynomials, we disregard $\infty$, which is always a rational fixed point.

\begin{theorem}
\label{123cycles}
Let $f(z)=z^2+c$ with $c \in \Q$.  Then
\begin{enumerate}
\item $f(z)$ has a rational point of period~$1$ (i.e., a rational fixed point) if and only if $c=1/4-\rho^2$ for some $\rho \in \Q$.  In this case, there are exactly two, $1/2+\rho$ and $1/2-\rho$, unless $\rho=0$, in which case they coincide.
\item $f(z)$ has a rational point of period~$2$ if and only if $c=-3/4-\sigma^2$ for some $\sigma \in \Q$, $\sigma \not=0$.  In this case, there are exactly two, $-1/2+\sigma$ and $-1/2-\sigma$ (and these form a 2-cycle).
\item $f(z)$ has a rational point of period~$3$ if and only if
	$$c = -\frac{\tau^6 + 2 \tau^5 + 4 \tau^4 + 8 \tau^3 + 9 \tau^2 + 4 \tau + 1}{4 \tau^2 (\tau+1)^2}$$
for some $\tau \in \Q$, $\tau \not= -1,0$.  In this case, there are exactly three,
\begin{align*}
	x_1 &=	\frac{\tau^3 + 2 \tau^2 + \tau + 1}{2 \tau (\tau+1)},	\\
	x_2 &=	\frac{\tau^3 - \tau - 1}{2 \tau (\tau+1)},	\\
	x_3 &=	-\frac{\tau^3 + 2 \tau^2 + 3 \tau + 1}{2 \tau (\tau+1)},
\end{align*}
and these are cyclically permuted by $f(z)$.
\end{enumerate}
\end{theorem}

\begin{proof}
This is a restatement of Theorems~1 and~3 in~\cite{walderusso}, except for the claim in~(3) that there can be at most one $3$-cycle, which is part of Theorem~3 in~\cite{morton3}.
\end{proof}

A criterion was given in~\cite{walderusso} for determining whether $z^2+c$ has rational points of period~$1$ {\em and} rational points of period~$2$.
The next theorem gives a more explicit criterion, and also shows that certain other combinations are impossible.

\begin{theorem}
\label{1and2}
Let $f(z)=z^2+c$ with $c \in \Q$.  Then
\begin{enumerate}
\item $f(z)$ has rational points of period~$1$ and rational points of period~$2$ if and only if
	$$c = -\frac{3 \mu^4 + 10 \mu^2 + 3}{4(\mu^2-1)^2}$$
for some $\mu \in \Q$, $\mu \not=-1,0,1$.
In this case the parameters $\rho$ and $\sigma$ of Theorem~\ref{123cycles} are
	$$\rho = - \frac{\mu^2+1}{\mu^2-1} \qquad \sigma = \frac{2 \mu}{\mu^2-1}.$$
\item If $f(z)$ has rational points of period $3$, it cannot have any rational points of period $1$ or $2$.
\end{enumerate}
\end{theorem}

\begin{proof}[Proof of Theorem~\ref{1and2}]
By Theorem~1, $z^2+c$ has rational points of period~$1$ and~$2$ if and only if
	$$c = 1/4 - \rho^2 = -3/4 - \sigma^2,$$
where $\rho,\sigma \in \Q$ with $\sigma \not=0$.
The curve in the $(\rho,\sigma)$-plane described by this equation is a conic with a rational point $(1,0)$, so it is birational to $\PP^1$ over $\Q$, with the rational function $\mu = (1-\rho)/\sigma$ giving the birational map.
Solving for $c$, $\rho$, and $\sigma$ in terms of $\mu$ gives the result.
The values of $\mu$ not allowed are $-1,1,0,\infty$, because these correspond to the two points at infinity on the conic and the two points where $\sigma=0$.

If $z^2+c$ has rational points of period $1$ and $3$, then
	$$c = 1/4 - \rho^2 = -\frac{\tau^6 + 2 \tau^5 + 4 \tau^4 + 8 \tau^3 + 9 \tau^2 + 4 \tau + 1}{4 \tau^2 (\tau+1)^2},$$
so $(\tau,2 \tau(\tau+1) \rho)$ is a point on the hyperelliptic curve
	$$C: \quad y^2 = x^6 + 2 x^5 + 5 x^4 + 10 x^3 + 10 x^2 + 4 x + 1.$$
By~\cite{elkies}, this is also an equation for the modular curve $X_1(18)$.
But $X_1(18)$ has only six rational points (these are all cusps), so the only rational points on $C$ besides the two points at infinity (on the nonsingular model) are $(-1,1)$, $(-1,-1)$, $(0,1)$, and $(0,-1)$.
These do not give rise to a valid pair $(\tau,\rho)$, since $\tau$ is not allowed to be~$0$ or~$-1$ in Theorem~\ref{123cycles}.
Hence it is impossible for there to exist rational points of period~$1$ and~$3$.

Similarly, if $z^2+c$ has rational points of period~$2$ and~$3$, then
	$$c = -3/4 - \sigma^2 = -\frac{\tau^6 + 2 \tau^5 + 4 \tau^4 + 8 \tau^3 + 9 \tau^2 + 4 \tau + 1}{4 \tau^2 (\tau+1)^2},$$
so $(\tau,2 \tau(\tau+1) \sigma)$ is a point on
	$$C': \quad y^2 = x^6 + 2 x^5 + x^4 + 2 x^3 + 6 x^2 + 4 x + 1.$$
This curve is $X_1(13)$, since if we the dehomogenize the model
	$$x_1^2 x_2^2 - x_1 x_2^3 - x_1 x_2 x_3^2 + x_1 x_3^3 + x_2^3 x_3 - x_2^2 x_3^2$$
of $X_1(13)$ in~\cite{billingmahler} by setting $x_3=1$, we find that the discriminant of the resulting quadratic in $x_1$ is
	$$x_2^6 + 2 x_2^5 + x_2^4 + 2 x_2^3 + 6 x_2^2 + 4 x_2 + 1.$$
The curve $X_1(13)$ also has exactly six rational points (again cusps), so the only rational points on $C'$ besides the two points at infinity are $(-1,1)$, $(-1,-1)$, $(0,1)$, and $(0,-1)$.
Again this implies that $z^2+c$ cannot have both rational points of period~$2$ and~$3$, since $\tau$ is not allowed to be $0$ or $-1$.
\end{proof}

The curves were originally recognized as $X_1(13)$ and $X_1(18)$ by computing enough invariants (such as the genus, automorphism group, primes of bad reduction, and Mordell-Weil group of the Jacobian) that the result could be guessed.

\section{Preperiodic Points}
\label{preperiodic}

If $m$ and $n$ are positive integers, then a {\em point of type $m_n$} for $f(z)$ is a preperiodic point that enters an $m$-cycle after $n$ iterations.
For example, $3/4$ is a point of type $3_2$ for $f(z)=z^2-29/16$, since its orbit is
	$$3/4, -5/4, -1/4, -7/4, 5/4, -1/4, -7/4, \ldots.$$

\begin{theorem}
\label{1_2}
Let $f(z)=z^2+c$ with $c \in \Q$.  Then
\begin{enumerate}
\item For each $m \ge 1$, $x$ is a rational point of type $m_1$ for $f(z)$ if and only if $-x$ is a nonzero rational point of period $m$.  The number of rational points of type $m_1$ equals the number of rational points of period $m$, except when  $c=-1$, $m=2$, or $c=0$, $m=1$, in which case there is one less.

\item $f(z)$ has rational points of type $1_2$ if and only if 
	$$c = \frac{-2(\eta^2+1)}{(\eta^2-1)^2}$$
for some $\eta \in \Q$, $\eta \not=-1,1$.
In this case, there are exactly $2$ such points, $\frac{2 \eta}{\eta^2-1}$ and $-\frac{2 \eta}{\eta^2-1}$, unless $\eta=0$ ($c=-2$), in which they coincide.
The parameter $\rho$ of Theorem~\ref{123cycles} is
	$$\rho = - \frac{\eta^2+3}{2(\eta^2-1)}.$$

\item $f(z)$ has rational points of type $2_2$ if and only if
	$$c = \frac{-\nu^4 - 2 \nu^3 - 2 \nu^2 + 2 \nu - 1}{(\nu^2-1)^2}$$
for some $\nu \in \Q$, $\nu \not=-1,0,1$.
In this case, there are exactly $2$ such points, $\frac{\nu^2+1}{\nu^2-1}$ and $-\frac{\nu^2+1}{\nu^2-1}$.
The parameter $\sigma$ of Theorem~\ref{123cycles} is
	$$\sigma = \frac{\nu^2 + 4 \nu -1}{2(\nu^2-1)}.$$

\item $f(z)$ has rational points of type $3_2$ if and only if $c=-29/16$.  For $c=-29/16$, the rational points of type $3_2$ are $3/4$ and $-3/4$.

\item If $f(z)$ has rational points of type $m_2$ with $1 \le m \le 3$, then there are no rational points of period $b \le 3$ unless $b=a$.

\item $f(z)$ cannot have rational points of type $1_n$, $2_n$, or $3_n$ for any $n \ge 3$.
\end{enumerate}
\end{theorem}

Together, Theorems~\ref{123cycles}, \ref{1and2}, and~\ref{1_2}, show that the subgraph of $\PrePer(z^2+c,\Q)$ consisting of finite rational preperiodic points that eventually end up in a cycle of length less than or equal to~$3$ is isomorphic to one of the graphs in Figure~1.
In particular, if Conjecture~\ref{periodissmall} holds, then these are the only possibilities for the full graph $\PrePer(z^2+c,\Q)$ (with $\infty$ deleted), so we get the following.

\begin{cor}
\label{preperiodiccor}
If Conjecture~\ref{periodissmall} holds, then $\#\PrePer(f,\Q) \le 9$ for all quadratic polynomials $f \in \Q[z]$.
\end{cor}
(The bound $9$ includes the fixed point at infinity.)  It should be remarked that many of the graphs in Figure~1 occur infinitely often: the only ones that are unique are the graphs for $c=-1$, $1/4$, $0$, $-2$, and $-29/16$.
Also, for each of the particular values of $c$ in the figure, it is easy to prove unconditionally that the graph shown is the full graph of finite rational preperiodic points.
In fact, this can be done even without explicit use of height functions: the rational preperiodic points are bounded in every absolute value, which means that they are multiples of $1/N$ for some fixed $N$ and are contained in some finite interval of $\R$, so that one need only evaluate $z^2+c$ at a finite number of points to compute the graph.
(In~\cite{callgoldstine}, similar arguments are used to bound the number of preperiodic points of $z^2+c$ in terms of the number of prime divisors of the denominator of $c$.)

\begin{proof}[Proof of Theorem~\ref{1_2}]
If $x$ is a rational point of type $m_1$, then $f(f^m(x))=f(x)$, but $f^m(x) \not= x$, so $f^m(x)=-x$.
Since every iterate of $x$ except $x$ itself is periodic of period $m$, $-x$ must be a rational point of period $m$, and it must be nonzero, since otherwise $x$ would be periodic as well.
Conversely, if $-x$ is a nonzero rational point of period $m$, then $f(x)=f(-x)$ is also a rational point of period $m$, but $f^m(x)=f^{m-1}(f(x))=f^{m-1}(f(-x))=-x \not=x$, so $x$ is a rational point of type $m_1$.

If $0$ is preperiodic for $z^2+c$, $c \in \Q$, then $c$ must be integral at each prime $p$, because otherwise the orbit of $0$ would diverge to infinity $p$-adically.
Thus $c \in \Z$.
Moreover, $|c| \le 2$, since otherwise the orbit of $0$ diverges to $+\infty$ with respect to the ordinary absolute value.
[This argument can be understood as follows: if $0$ is preperiodic, then its orbit is bounded in every absolute value.  (In fact, the converse holds as well.)  In other words, $0$ is in the Mandelbrot set and also in the $p$-adic Mandelbrot set (which is the closed unit ball in $\Qbar_p$) for each prime $p$.]
Checking $c=-2,-1,0,1,2$, we find that $0$ is preperiodic only for $c=-2,-1,0$, and periodic only for $c=-1$ and $c=0$, with period $2$ and $1$, respectively.
Thus the number of rational points of type $m_1$ equals the number of rational points of period~$m$, except for $c=-1$, $m=2$, and $c=0$, $m=1$, when it is one less.
This completes the proof of~(1).

\medskip

If $z^2+c$ has a rational point $r$ of type $1_2$, then $c=1/4-\rho^2$ by Theorem~\ref{123cycles}, and $r^2+c$ is the negative of one of the fixed points, since these are exactly the points of type $1_1$, provided that the fixed point is not $0$.
Without loss of generality, we may assume $r^2+c=-(1/2+\rho)$ with $\rho \not=-1/2$, so
	$$r^2 = -(1/2+\rho)-(1/4-\rho^2) = \rho^2 - \rho - 3/4.$$
This conic in the $(r,\rho)$-plane has a rational point $(0,-1/2)$, so the rational function $\eta=r/(1/2+\rho)$ defines a birational map from it to $\PP^1$.
Substituting $r=(1/2+\rho)\eta$ in the equation of the conic, we solve and find
	$$\rho=- \frac{\eta^2+3}{2(\eta^2-1)}, \quad r=-\frac{2 \eta}{\eta^2-1}, \quad c= \frac{-2(\eta^2+1)}{(\eta^2-1)^2}.$$
The values of $\eta$ not allowed are $1,-1,\infty$, because these correspond to the two points at infinity on the conic, and the point where $\rho=-1/2$.
But for the good values of $\eta$, $\frac{2 \eta}{\eta^2-1}$ is another rational point of type $1_2$, unless $\eta=0$, in which case they coincide.

If there were yet another point $s$ of type $1_2$, then $s^2+c$ would have to be the negative of the {\em other} fixed point: $s^2+c=-(1/2-\rho)$.  Then
	$$s^2 = \frac{-\eta^4+2 \eta^2 + 3}{(\eta^2-1)^2}.$$
Letting $t=(\eta^2-1)s$, we find that we have a point on the genus $1$ curve
	$$t^2=-\eta^4+2 \eta^2 + 3.$$
Since the curve has a rational point $(\eta,t)=(1,2)$, it is birational to an elliptic curve.
In fact, the following are inverse rational maps between the curve and the elliptic curve $E_{24}:y^2=x^3-x^2+x$ in Weierstrass form:
\begin{align*}
	(\eta,t)& \mapsto \left( \frac{t+2}{(\eta-1)^2} , \frac{2 \eta^3 - 2 \eta^2 - 2\eta - 6 - 4t}{2(\eta-1)^3} \right),	\\
	(x,y)	& \mapsto \left( \frac{x^2-2y-1}{x^2+1} , \frac{2 x^4 - 4 x^3 + 8 xy + 4x -2}{(x^2+1)^2} \right).
\end{align*}
The elliptic curve $E_{24}$ is curve 24A4 in Cremona's tables~\cite{cremona} and 24A in Antwerp IV~\cite{antwerp4}, and
	$$E_{24}(\Q) = \{\calO,(0,0),(1,1),(-1,1)\}.$$
These correspond to the points $(\eta,t)=(1,2),(-1,-2),(-1,2),(1,-2)$, respectively.
But $\eta=1$ and $\eta=-1$ are forbidden, so there can be no further points of type $1_2$.

\medskip

If $z^2+c$ has a rational point $r$ of type $2_2$, then $c=-3/4-\sigma^2$ for some nonzero $\sigma \in \Q$ by Theorem~\ref{123cycles}, and $r^2+c$ is the negative of a nonzero point of period~$2$.
Without loss of generality, we may assume $r^2+c=-(-1/2+\sigma)$ with $\sigma \not=0,1/2$, so
	$$r^2 = -(-1/2+\sigma)-(-3/4-\sigma^2) = \sigma^2 - \sigma + 5/4.$$
This conic in the $(r,\sigma)$-plane has a rational point $(1,1/2)$, so the rational function $\nu=(r-1)/(1/2-\sigma)$ defines a birational map from it to $\PP^1$.
Substituting $r=1+(1/2-\sigma)\nu$ in the equation of the conic, we solve and find
	$$\sigma = \frac{\nu^2 +  4\nu -1}{2(\nu^2-1)}, \quad r=-\frac{\nu^2+1}{\nu^2-1}, \quad c= \frac{-\nu^4 - 2 \nu^3 - 2 \nu^2 + 2 \nu - 1}{(\nu^2-1)^2}.$$
The values of $\nu$ not allowed are $1,-1,0,\infty$, because these correspond to the two points at infinity on the conic and the two points where $\sigma=1/2$.  (Actually $-2+\sqrt{5}$ and $-2-\sqrt{5}$ are not allowed either since they correspond to the points where $\sigma=0$, but these are not rational, so we need not consider these.)
For the good values of $\nu$, $\frac{\nu^2+1}{\nu^2-1}$ is another rational point of type $2_2$.

If there were yet another point $s$ of type $2_2$, then $s^2+c$ would have to be the negative of the {\em other} point of period~$2$: $s^2+c=-(-1/2-\sigma)$.  Then
	$$s^2 = \frac{2(\nu^4+2 \nu^3-2 \nu+1)}{(\nu^2-1)^2}.$$
Letting $t=(\nu^2-1)s$, we find that we have a point on the genus $1$ curve
	$$t^2=2(\nu^4+2 \nu^3-2 \nu+1).$$
Since this curve has a rational point $(\nu,t)=(1,2)$, it is birational to an elliptic curve.
The following are inverse rational maps between the curve and the elliptic curve $E_{40}:y^2=x^3-2x+1$:
\begin{align*}
	(\nu,t)	& \mapsto \left( \frac{t+2 \nu^2}{(\nu-1)^2} , -\frac{3 \nu^3 + 2\nu t + 3 \nu^2 -3 \nu+1}{(\nu-1)^3} \right),	\\
	(x,y)	& \mapsto \left( \frac{x^2-2y}{x^2-4x+2} , \frac{2x^4-8x^2 y+8x^3+8xy-24x^2+24x-8}{(x^2-4x+2)^2} \right).
\end{align*}
The elliptic curve $E_{40}$ is curve 40A3 in~\cite{cremona} (40A in~\cite{antwerp4}), and
	$$E_{40}(\Q) = \{\calO,(0,1),(1,0),(0,-1)\}.$$
These correspond to the points $(\nu,t)=(1,2),(-1,-2),(-1,2),(1,-2)$, respectively.
But $\nu=1$ and $\nu=-1$ are forbidden, so there can be no further points of type $2_2$.

\medskip

We will postpone the discussion of points of type $3_2$ until the next section.
The case $a=3$ of part~(5) will follow from this, so for now we prove part~(5) for $a=1$ and $a=2$.
Note that in each of these cases, there can be no points of order~$3$, by Theorem~\ref{1and2}.

If $f(z)$ has a rational point $r$ of type $1_2$ then $c=1/4-\rho^2$ and $r^2+c$ is the negative of a nonzero fixed point, as before: without loss of generality $r^2+c=-(1/2-\rho)$.  If $f(z)$ also has a rational point of period~$2$, then
	$$c = -\frac{3 \mu^4 + 10 \mu^2 + 3}{4(\mu^2-1)^2}, \quad \rho = - \frac{\mu^2+1}{\mu^2-1},$$
for some $\mu \in \Q$, $\mu \not=-1,0,1$, and we get
	$$r^2 = \frac{5 \mu^4 + 14 \mu^2 -3}{4(\mu^2-1)^2}.$$
Letting $t=2(\mu^2-1)r$, we obtain a rational point on the genus $1$ curve
	$$t^2=5 \mu^4 + 14 \mu^2 -3.$$
Since this curve has a rational point $(\mu,t)=(1,4)$, it is birational to an elliptic curve.
The following are inverse rational maps between the curve and the elliptic curve $E_{15}:y^2+xy+y=x^3+x^2$:
\begin{align*}
	(\mu,t)	& \mapsto \left( \frac{t+\mu^2+4\mu-1}{2(\mu-1)^2} , -\frac{2 \mu^3 + \mu t + \mu^2 +2 \mu - 1}{(\mu-1)^3} \right),	\\
	(x,y)	& \mapsto \left( \frac{x^2-2y+x-1}{x^2-x-1} , \frac{4x^4-12x^2 y+8x^3-8xy+8x^2+4x-8y-4}{(x^2-x-1)^2} \right).
\end{align*}
The elliptic curve $E_{15}$ is curve 15A8 in~\cite{cremona} (15A in~\cite{antwerp4}), and
	$$E_{15}(\Q) = \{\calO,(0,0),(-1,0),(0,-1)\}.$$
These correspond to the points $(\mu,t)=(1,4),(1,-4),(-1,4),(-1,4)$, respectively.
But $\mu=1$ and $\mu=-1$ are forbidden, so it is impossible for $f(z)$ to have a rational point of type $1_2$ and a rational point of period $2$.

Similarly, if $f(z)$ has a rational point $r$ of type $2_2$, then $c=-3/4-\sigma^2$ and $r^2+c$ is the negative of a nonzero period~$2$ point, as before: without loss of generality $r^2+c=-(-1/2+\sigma)$.  If $f(z)$ also has a rational point of period~$1$, then
	$$c = -\frac{3 \mu^4 + 10 \mu^2 + 3}{4(\mu^2-1)^2}, \quad \sigma = \frac{2\mu}{\mu^2-1},$$
for some $\mu \in \Q$, $\mu \not=-1,0,1$, and we get
	$$r^2 = \frac{5 \mu^4 - 8 \mu^3 + 6 \mu^2 +8 \mu + 5}{4(\mu^2-1)^2}.$$
Letting $t=2(\mu^2-1)r$, we obtain a rational point on the genus $1$ curve
	$$t^2=5 \mu^4 - 8 \mu^3 + 6 \mu^2 +8 \mu + 5.$$
Since this curve has a rational point $(\mu,t)=(1,4)$, it is birational to an elliptic curve.
The following are inverse rational maps between the curve and the elliptic curve $E_{17}:y^2+xy+y=x^3-x^2-x$:
\begin{align*}
	(\mu,t)	& \mapsto \left( \frac{t+\mu^2+3}{2(\mu-1)^2} , -\frac{3 \mu^3 + \mu t - 5 \mu^2 + 9 \mu + 1}{2(\mu-1)^3} \right),	\\
	(x,y)	& \mapsto \left( \frac{x^2-2y-x-1}{x^2-x-1} , \frac{4x^4-4x^2 y-4x^3-8xy-4x-4}{(x^2-x-1)^2} \right).
\end{align*}
The elliptic curve $E_{17}$ is curve 17A4 in~\cite{cremona} (17A in~\cite{antwerp4}), and
	$$E_{17}(\Q) = \{\calO,(0,0),(1,-1),(0,-1)\}.$$
These correspond to the points $(\mu,t)=(1,4),(1,-4),(-1,4),(-1,-4)$, respectively.
But $\mu=1$ and $\mu=-1$ are forbidden, so it is impossible for $f(z)$ to have a rational point of type $2_2$ and a rational point of period~$1$.
This completes the proof of part~(5) for $a \le 2$.

\medskip

It follows from part~(4) and the calculation of preperiodic points for $z^2-29/16$ (which we have postponed until the next section), that $f(z)$ cannot have rational points of type $3_n$ for $n \ge 3$.
Therefore we now prove that $f(z)$ cannot have rational points of type $1_n$ or $2_n$ for $n \ge 3$.
It suffices to consider $n=3$, since iterating $f$ on a rational point with higher $n$ eventually yields a point with $n=3$.
If $f(z)$ has a rational point $q$ of type $1_3$, then $q^2+c$ is a rational point $r$ of type $1_2$, and by part~(2), we have
	$$c = \frac{-2(\eta^2+1)}{(\eta^2-1)^2}, \quad r=-\frac{2 \eta}{\eta^2-1},$$
where $\eta \not= -1,1$, so
	$$q^2 = \frac{2(\eta^3+\eta^2-\eta+1)}{(\eta^2-1)^2}.$$
Letting $t=2(\eta^2-1)q$, we obtain a rational point on the elliptic curve
	$$t^2=2(\eta^3+\eta^2-\eta+1).$$
The linear change of coordinates $\eta=2x-1$, $t=4y+2$ transforms this into
	$$E_{11}: y^2 + y = x^3 - x^2,$$
which is curve 11A3 in~\cite{cremona} (11A in~\cite{antwerp4}), and which is also the modular curve $X_1(11)$.
This curve has five rational points
	$$\calO,(0,0),(0,-1),(1,0),(1,-1)$$
the last four of which correspond to finite points
	$$(\eta,t)=(-1,2),(-1,-2),(1,2),(1,-2)$$
but $\eta=1$ and $\eta=-1$ are forbidden, so $f(z)$ cannot have rational points of type $1_3$.

If $f(z)$ has a rational point $q$ of type $2_3$, then $q^2+c$ is a rational point $r$ of type $2_2$, and by part~(2), we have
	$$c = \frac{-\nu^4 - 2 \nu^3 - 2 \nu^2 + 2 \nu - 1}{(\nu^2-1)^2}, \quad r = \frac{\nu^2+1}{\nu^2-1},$$
for some $\nu \in \Q$, $\nu \not=-1,0,1$,
so
	$$q^2 = \frac{2(\nu^3+\nu^2-\nu+1)}{(\nu^2-1)^2}.$$
This is the same equation found in the last paragraph (with $\nu$ instead of $\eta$).
Again $\nu=1$ and $\nu=-1$ are forbidden, so $f(z)$ cannot have rational points of type $2_3$.
\end{proof}

\section{Rational points on the genus~$2$ curve $C_1(3_2)$}
\label{3_2}

In this section, we complete the proof of Theorem~\ref{1_2} by proving part~(4).
First we find a nice model for the curve $C_1(3_2)$ which classifies $z^2+c$ together with a point of type $3_2$
If $z^2+c$ has a rational point $r$ of type $3_2$, then $r^2+c$ is the negative of a nonzero rational point of period~$3$, so by Theorem~\ref{123cycles} we have, without loss of generality,
\begin{align*}
	c &= -\frac{\tau^6 + 2 \tau^5 + 4 \tau^4 + 8 \tau^3 + 9 \tau^2 + 4 \tau + 1}{4 \tau^2 (\tau+1)^2}	\\
	r^2+c &= -\frac{\tau^3 + 2 \tau^2 + \tau + 1}{2 \tau (\tau+1)},
\end{align*}
for some $\tau \in \Q$, $\tau \not= -1,0$.
These imply that
	$$r^2 = \frac{\tau^6 - 2 \tau^4 + 2 \tau^3 + 5 \tau^2 + 2 \tau + 1}{4 \tau^2 (\tau+1)^2},$$
so $(\tau,2 \tau(\tau+1) r)$ is a rational point on the curve
	$$\calC: \quad y^2 = g(x)$$
where
	$$g(x) \defequal x^6 - 2 x^4 + 2 x^3 + 5 x^2 + 2 x + 1.$$

Since $g(x)$ has no rational zeros, $\calC$ cannot be put in the form $y^2=\text{(quintic in $x$)}$ over $\Q$.
Furthermore it can be shown using the same methods as in~\cite{fivecycle} that the Jacobian $J$ of $\calC$ is absolutely simple with $\End J=\Z$, so that $J$ is not a quotient of the Jacobian of any modular curve over $\C$.
Hence we will apply the general method of~\cite{fivecycle} to find the rational points.
(We will be brief in this section; the reader is encouraged to consult~\cite{fivecycle} for more details.)

Note that $\calC$ has six obvious affine points: $Q^+=(-1,1)$, $Q^-=(-1,-1)$, $R^+=(0,1)$, $R^-=(0,-1)$, $S^+=(1,3)$, and $S^-=(1,-3)$.
Since $\deg g$ is even, $\calC$ has two points at infinity (on the nonsingular model), and these are rational since the leading coefficient of $g(x)$ is a square.
The rational function $y/x^3$ takes values~1 and~$-1$ at these two points, which we call $\infty^+$ and $\infty^-$, respectively.
We will eventually show that these eight points are the only rational points on $\calC$.
First we compute the structure of the Mordell-Weil group $J(\Q)$.
Elements of $J(\Q)$ will be written either as $[D]$ where $D$ is a degree zero divisor, or as $[P_1+P_2]$ with $P_1,P_2 \in \calC$, which is (in abuse of notation) identified with $[P_1+P_2-\infty^+-\infty^-]$.

\begin{table}
\begin{center}
\begin{tabular}{|c|c|c|}
Element		& Definition				& Norm \\ \hline \hline
$u_1$		& $(T^4-T^3-T^2+2T+1)/2$		& $1$	\\ \hline
$u_2$		& $(T^4-T^3-T^2+4T+1)/2$		& $1$	\\ \hline
$-1$		& $-1$					& $1$	\\ \hline
$\alpha$	& $(T^5-2T^3+T^2+7T+3)/2$		& $2^3$	\\ \hline
$\beta_1$	& $(T^5-5T^3+5T^2+6T-2)/2$		& $743$ \\ \hline
$\beta_2$	& $(T^5+8T^4-10T^3-3T^2+35T+13)/2$	& $743^2$ \\ \hline
$\beta_3$	& $(-10T^5+9T^4+14T^3-33T^2-21T+18)/2$	& $743^2$ \\ \hline
\end{tabular}
\end{center}
\caption{Some elements of $L$.}
\label{elements}
\end{table}

\begin{prop}
\label{rank}
$J(\Q) \isom \Z$.
\end{prop}

\begin{proof}
The discriminant of $g(x)$ is $-2^{12} \cdot 743$, so the set of bad primes of $J$ is contained in $S \defequal \{2,743,\infty\}$.
(In fact, $\calC$ and $J$ have good reduction at~$2$, because substituting $y=2z+x^3+x+1$ and dividing by~$4$ yields the model
	$$z^2 + z x^3 + z x + z + x^4 = x^2,$$
which has bad reduction only at $743$.)
We compute $\# J(\F_3) = 27$ and $\# J(\F_5) = 43$, so $J(\Q)\tors$ is trivial.
Thus the divisor class $[\infty^+ - \infty^-]$ has infinite order, and $J(\Q)$ has rank at least~$1$.

Let $L$ be the field $\Q[T]/(g(T))$, and for each prime $p$ including $\infty$ let $L_p=\Q_p[T]/(g(T))$.
PARI tells us that the splitting field of $g(x)$ over $\Q$ has Galois group $S_6$, that the class number of $L$ is~$1$, and that the unit group is of rank~$2$, generated by $u_1$, $u_2$, and $-1$, where $u_1$ are $u_2$ are defined in Table~\ref{elements}.
Furthermore the ramified primes~$2$ and $743$ factor in $L$ as $-\alpha^2 u_1$ and $\beta_1^2 \beta_2 \beta_3$, where $\alpha$, $\beta_1$, $\beta_2$, $\beta_3$ are irreducible elements also defined in Table~\ref{elements}.
We have $L_{743} \isom E \times F \times F$, where $E$ and $F$ are quadratic extensions of $\Q_{743}$ with $E$ ramified and $F$ unramified.

As in~\cite{fivecycle}, there is a map
	$$(x-T): J(\Q) \rightarrow L^\ast/L^{\ast 2}\Q^\ast.$$
The index of $2J(\Q)$ in the kernel is~$2$ by~\cite[Proposition~5]{fivecycle}, since $g(x)$ has Galois group $S_6$.
So to prove that $J(\Q)$ has rank~$1$, it suffices to prove that the image is trivial.
The image is contained in the subgroup $H$ of elements in the kernel of the norm $L^\ast/L^{\ast 2}\Q^\ast \rightarrow \Q^\ast/\Q^{\ast 2}$ which are ``unramified outside $S$,'' by~\cite[Proposition~3]{fivecycle}.
In fact, it is contained in the subgroup $H_0$ of $H$ consisting of elements which map in $L_{743}^\ast/L_{743}^{\ast 2} \Q_{743}$ into the image of
	$$(x-T): J(\Q_{743}) \rightarrow L_{743}^\ast/L_{743}^{\ast 2}\Q_{743}^\ast.$$
(It will turn out that we will not need the information from the other primes in $S$.)

In our case, an $\F_2$-basis for the subgroup of $L^\ast/L^{\ast 2}$ unramified outside $S$ is the set of seven elements in Table~\ref{elements}.
The intersection with the kernel of the norm to $\Q^\ast/\Q^{\ast 2}$ is spanned by $u_1,u_2,-1,\beta_2,\beta_3$.
Since $u_1 = -2 \alpha^{-2} \in L^{\ast 2} \Q^\ast$, and similarly $\beta_2 \equiv \beta_3 \pmod{L^{\ast 2} \Q^\ast}$, we find that $H$ is spanned by (the images of) $u_2$ and $\beta_2$.

Let us now calculate the image of
	$$(x-T): J(\Q_{743}) \rightarrow L_{743}^\ast/L_{743}^{\ast 2}\Q_{743}^\ast.$$
The~$15$ nonzero elements of the $2$-torsion subgroup $J[2]$ are of the form $[P_1+P_2]$ where $P_1$ and $P_2$ are distinct points on $\calC$ of the form $(r,0)$ with $g(r)=0$.
Since $g(x)$ factors into three quadratics over $\Q_{743}$, such $[P_1+P_2]$ will be stable under $\Gal(\Qbar_{743}/\Q_{743})$ if and only if the two corresponding zeros of $g(x)$ are zeros of the same quadratic factor.
Hence $\# J(\Q_{743})[2] = 3+1 = 4.$
Multiplication-by-2 is locally Haar measure-preserving on $J(\Q_{743})$ so
	$$\# J(\Q_{743})/2J(\Q_{743}) = \# J(\Q_{743})[2] = 4.$$
Finally, $g(x)$ has no zero in $\Q_{743}$, and no partition of the six zeros of $g(x)$ into two 3-element subsets can be stable under $\Gal(\Qbar_{743}/\Q_{743})$ because of the decomposition of $L_{743}$, so by~\cite[Proposition~5]{fivecycle}, $2J(\Q_{743})$ has index~$2$ in $\ker(x-T)$, and
\begin{equation}
\label{image}
	\# (x-T)(J(\Q_{743})) = \# J(\Q_{743})/\ker(x-T) = 2.
\end{equation}
Since the Legendre symbol $\left(\frac{33}{743}\right)$ is~$1$, Hensel's Lemma implies that $(2,\sqrt{33}) \in \calC(\Q_{743})$.  (Fix a square root.)
We claim that $[(2,\sqrt{33}) - \infty^-]$ generates $J(\Q_{743})/\ker(x-T)$, which is equivalent to $2-T$ generating $(x-T)(J(\Q_{743}))$.
Because of~(\ref{image}), it suffices to show $2-T \not\in L_{743}^{\ast 2} \Q_{743}^\ast$.

In fact, we will show more: that $2-T$, $u_2$ and $\beta_2$ are $\F_2$-independent in $L_{743}^\ast/L_{743}^{\ast 2}\Q_{743}^\ast$.
If not, then since $\Q_{743}^\ast/\Q_{743}^{\ast 2}$ is generated by $-1$ and $743$, some product of $-1$, $743$, $2-T$, $u_2$, and $\beta_2$ involving at least one of the last three would be in $L_{743}^{\ast 2}$.
Recall that $L_{743} \isom E \times F \times F$, where $F=\Q_{743}(i)$ is the unramified quadratic extension of $\Q_{743}$, $i^2=-1$.
Since $-1$ is a square in $F$, we would find that the product of some nonempty subset of $\{743, 2-T, u_2, \beta_2\}$ would map to the trivial element of $F^\ast/F^{\ast 2} \times F^\ast/F^{\ast 2}$.
The condition that the $\beta_3$-adic valuation of the product be even forces $743$ not to be part of the product.
The condition that the $\beta_2$-adic valuation of the product be even forces $\beta_2$ not to be part of the product.
So some $\F_2$-combination of $2-T$ and $u_2$ is a square in both $F$'s.
By factoring $g(x)$ modulo $743$, we find that $T$ maps in the first $F$ to something congruent to $330+2i$ in the residue field $\F_{743^2}$, and in the second $F$ to something congruent to $458+44i$ in the residue field.
Checking the Legendre symbols at the norm from $\F_{743^2}$ to $\F_{743}$ of the image of $2-T$ and $u_2$ for each $F$, we find by Hensel's Lemma that $2-T$ is a square in the first $F$ but not the second, and that $u_2$ is not a square in either $F$.
This contradicts the existence of a relation between $2-T$ and $u_2$ in $F^\ast/F^{\ast 2} \times F^\ast/F^{\ast 2}$, and hence proves that $2-T$, $u_2$ and $\beta_2$ are $\F_2$-independent in $L_{743}^\ast/L_{743}^{\ast 2}\Q_{743}^\ast$.

But this means that $2-T$ generates the image of $J(\Q_{743})$ under $(x-T)$, and also that no element of $H$ can map in $L_{743}^\ast/L_{743}^{\ast 2}\Q_{743}^\ast$ into this image.
Thus $J(\Q)/\ker(x-T)$ is trivial, and the $\F_2$-dimension of $J(\Q)/2J(\Q)$ is at most~$1$.
We have already shown that $J(\Q)$ is torsion-free and of rank at least~$1$, so $J(\Q) \isom \Z$.
\end{proof}

Since the rank of $J(\Q)$ is~$1$, we can now apply the method of Chabauty and Coleman to bound the number of rational points on $\calC$.
We will work modulo powers of~$3$.
Over $\F_3$, $\calC$ has five affine points, so $\# \calC(\F_3) = 7$.
The mod~$3$ reductions of the eight known rational points $Q^+$, $Q^-$, $R^+$, $R^-$, $S^+$, $S^-$, $\infty^+$, $\infty^-$ are distinct, except that $R^+$ and $R^-$ both reduce to the Weierstrass point $(1,0) \in \calC(\F_3)$.
To show these are all, it suffices to prove the following:

\begin{prop}
Each point in $\calC(\F_3)$ is the mod~$3$ reduction of exactly one point in $\calC(\Q)$ except $(1,0)$ which is the reduction of exactly two.
\end{prop}

\begin{proof}
Let $D=[\infty^+-\infty^-]=[\infty^++\infty^+] \in J(\Q)$ (recall our abuse of notation), and let $\tildeD \in J(\F_3)$ be the mod~$3$ reduction of $D$.
Using the group law presented in~\cite{flynngrouplaw} (or alternatively by intersecting $\calC$ with carefully chosen cubics to obtain relations between divisor classes), we find that $9D=[Q^-+R^+]$ and $27D=[S^-+S^-]$.
Since $\#J(\F_3)=27$ and $9\tildeD \not= \calO \in J(\F_3)$, $J(\F_3)$ is cyclic of order $27$ and $\tildeD$ is a generator.
Let $\calM_3 \subset J(\Q_3)$ denote the kernel of reduction.
Then $D'=27 \cdot D=[S^-+S^-]$ is in $\calM_3$.
We do not know whether $D$ is a generator of $J(\Q)$, but if $E$ is a generator and $D=n \cdot E$, then $n$ is not divisible by $3$, since $\tildeD \not\in 3J(\F_3)$.
This is enough to imply that every $D_0 \in J(\Q) \cap \calM_3$ is (uniquely) expressible as $n \cdot D'$, where $n$ is a {\em $3$-adic} integer.

%
Suppose $P \in \calC(\Q)$ has the same mod~$3$ reduction as $Q^+$.
Then $[P+P]$ equals $[Q^++Q^+] + n \cdot D'$ for some $n \in \Z_3$.
Using the notation and methods of~\cite{fivecycle}, we find that the formal logarithm of $D'$ satisfies
	$$\ell_1 \equiv 3 \pmod{3^4}, \quad \ell_2 \equiv 75 \pmod{3^4},$$
and we compute a power series $\theta(n) \in \Z_3[[n]]$ which vanishes whenever $[Q^++Q^+] + n \cdot D'$ is of the form $[P+P]$.
In our case, we find $\theta(n) \equiv 3n \pmod{3^2}$.
Strassman's theorem~\cite[p. 62]{cassels} implies that such a power series can have at most one zero in $\Z_3$.
We already know that $n=0$ is a zero, so it the only one.
Thus $[P+P]=[Q^++Q^+]$, and since $J(\Q)$ has no 2-torsion, $P=Q^+$.
Applying the hyperelliptic involution, we deduce also that any $P \in \calC(\Q)$ with the same mod~$3$ reduction as $Q^-$ must equal $Q^-$.

Applying the same argument with $R^+$ instead of $Q^+$, we find $\theta(n) \equiv 6n \pmod{3^2}$, so we deduce similarly that any $P \in \calC(\Q)$ with the same mod~$3$ reduction as $R^+$ or $R^-$ must equal $R^+$ or $R^-$.
For $\infty^+$ instead of $Q^+$, we obtain $\theta(n) \equiv 6n \pmod{3^2}$, so any $P \in \calC(\Q)$ with the same mod~$3$ reduction as $\infty^+$ or $\infty^-$ must equal $\infty^+$ or $\infty^-$.

It remains to be shown that if $P \in \calC(\Q)$ reduces to $(1,0)$ modulo~$3$, then $P$ equals $S^+$ or $S^-$.
Although in theory the same technique as before could be used to bound the number of such rational points, it turns out that in this case even the terms up to degree~$7$ of the formal logarithm and exponential do not give sufficient $3$-adic precision.
(In other words, we cannot distinguish $\theta(n)$ from $0$ without using a larger number of terms.)
Qualitatively, the difficulty seems to be that $(1,0)$ is a Weierstrass point of $\calC$ over $\F_3$, so that the image of $\calC$ in $J$ under $P \mapsto [P+P]$ has six branches passing the origin $\calO$, one for each Weierstrass point, and we can therefore expect $\theta(n)$ to have a zero of multiplicity~$6$ and to be divisible by many factors of~$3$.
We will circumvent the need for more precision by substituting the embedding $P \mapsto [P+S^+]$ of $\calC$ into $J$ for the embedding $P \mapsto [P+P]$.

The rational function $t=y+3$ is a uniformizing parameter at $S^-=(1,-3)$, so there is a unique power series $\xi(t) \in \Q[[t]]$ which starts
	$$\xi(t) = 1 - 3t/8 - 31t^2/512 + 105t^3/16384 + 15269t^4/2097152+O(t^5)$$
such that $(\xi(t),-3+t)$ is a point on $\calC$ with coordinates in $\Q[[t]]$.
Since $t$ is also a uniformizing parameter on the the curve reduced modulo $3$, $\xi(t) \in \Z_3[[t]]$, and hence for any $t \in 3\Z_3$, the power series converges to give a point $P_t$ in $\calC(\Q_3)$.
This analytic parameterization gives all points in $\calC(\Q_3)$ which have the same mod~$3$ reduction as $S^-$.
Let $D_t=[P_t+S^+] \in \calM_3 \subset J(\Q_3)$.
The local parameters $s_1,s_2$ of~\cite{fivecycle} at $D_t$ are again given by power series in $\Z_3[[t]]$:
\begin{align*}
	s_1(t) &= t/16 + 9 t^2/1024 - 111 t^3/32768 - 8979 t^4/4194304 + O(t^5),\\
	s_2(t) &= t/16 + 21 t^2/1024 + 141 t^3/32768 + 1929 t^4/4194304 + O(t^5).
\end{align*}
Using the formulas in~\cite{fivecycle}, we find that the power series giving the formal logarithm at $D_t$ for $t=3n$ satisfy
\begin{align*}
	L_1(n) &\equiv 66 n + 54 n^3 \pmod{3^4}, \\
	L_2(n) &\equiv 66 n + 27 n^2 + 72 n^3 \pmod{3^4}.
\end{align*}
If $t=3n \in 3\Z_3$ is such that $P_t \in \calC(\Q)$, then $D_t \in J(\Q) \cap \calM_3$, so its logarithm $(L_1(n),L_2(n))$ must be a $3$-adic integer multiple of the logarithm $(\ell_1,\ell_2)$ of $D'$.
In particular the determinant
	$$\Delta(n) \defequal \begin{vmatrix} L_1(n) & L_2(n) \\ \ell_1 & \ell_2 \end{vmatrix}$$
must vanish whenever $P_t \in \calC(\Q)$.
The power series $\Delta(n)$ satisfies
	$$\Delta(n) \equiv 54 n + 27 n^3 \pmod{3^4},$$
so by Strassman's Theorem, $\Delta(n)=0$ for at most three values of $n \in \Z_3$.
We already know that $\Delta(0)=\Delta(2)=0$, since $P_0=(1,-3)$ and $P_6=(1,3)$.
At $n=1$, $P_3$ is a Weierstrass point $W$ defined over $\Q_3$ and not $\Q$, but $\Delta(1)=0$ nevertheless: $2D_3=[S^++S^+]+[W+W]=[S^++S^+]=-D'$, so the logarithm of $D_3$ is a multiple of the logarithm of $D'$.
Thus $\Delta(n) \not=0$ for $n \not=0,1,2$, and it follows that any $P \in \calC(\Q)$ which reduces to $(1,0) \in \calC(\F_3)$ must be one of $P_0=S^-$, $P_1=W$, $P_2=S^+$.
Of these, only $S^+$ and $S^-$ are actually rational, so we are done.
\end{proof}

Thus there are only the eight points on $\calC$.
If $z^2+c$ has a rational point of type $3_2$ the corresponding value of $\tau=x$ must be one of $-1,0,1,\infty$.
But $\tau=-1,0,\infty$ do not correspond to a valid polynomial $z^2+c$ in Theorem~\ref{123cycles}, so the only possibility is $\tau=1$, which makes $c=-29/16$.
Computing the full graph $\PrePer(f,\Q)$ for $f(z)=z^2-29/16$ (which can be done as in the remarks following Corollary~\ref{preperiodiccor}) completes the proof of part~(4) of Theorem~\ref{1_2}.

\section*{Acknowledgements}

I thank Noam Elkies for computing an equation of $X_1(18)$ for me, and E.\ V.\ Flynn for making available his formulas for Jacobians of genus~$2$ curves.


\end{document}